\documentclass[a4paper]{article}
\usepackage{amsmath,amsthm,amsfonts,amssymb,mathrsfs}
\usepackage[latin1]{inputenc}

\usepackage{hyperref}

\newcommand{\N}{\mathbb{N}}

\newcommand{\R}{\mathbb{R}}

\newcommand{\T}{\mathbb{T}}
\newcommand{\Ex}{\mathbb{E}}
\newcommand{\Prob}{\mathbb{P}}

\newtheorem{theorem}{Theorem}

\theoremstyle{remark}
\newtheorem*{remark}{Remark}
\newtheorem*{example}{Example}

\begin{document}

\title{Supercritical general branching processes conditioned on extinction are subcritical}

\author{Peter Jagers\footnotemark[1] \footnotemark[3] \and Andreas~N.~Lager{\aa}s\footnotemark[2] \footnotemark[4]}

\footnotetext[1]{Supported by the Swedish Research Council}
\footnotetext[2]{Supported by the Science Faculty of Gothenburg University through the Platform for Theoretical Biology}
\footnotetext[3]{E-mail: jagers@chalmers.se}
\footnotetext[4]{E-mail: norand@chalmers.se}

\maketitle

\begin{abstract}
It is well known that a simple, supercritical
Bienaym\'e-Galton-Watson process turns into a subcritical such
process, if conditioned to die out. We prove that the corresponding
holds true for general, multi-type branching, where child-bearing
may occur at different ages, life span may depend upon reproduction,
and the whole course of events is thus affected by conditioning upon
extinction.\\

\end{abstract}

\section{Introduction}

The theory of branching processes was born out of Galton's famous
family extinction problem. Later, interest turned to populations not
dying out and their growth and stabilisation. In more recent years,
extinction has retaken a place in the foreground, for reasons from
both conservation and evolutionary biology. The time and path to
extinction of subcritical general populations was studied in
\cite{pnas}. Here, time structure is crucial, and life spans and
varying bearing ages cannot be condensed into simple, generation
counting Bienaym\'e-Galton-Watson processes. Thus, the question
arises whether (non-critical)  general branching populations bound
for extinction must behave like subcritical populations.

We answer this in the affirmative: a general, multitype branching
process conditioned to die out, remains a branching process, but one
almost surely dying out. If the original process was supercritical
but with a strictly positive risk of extinction, the new process is
subcritical.

Formulated in such a loose manner, this fact belongs to the folklore
of branching, but actually it has been  proved only for
Bienaym\'e-Galton-Watson processes, \cite{AthreyaNey:1972}, p. 52. A
moment's afterthought tells us that it remains true for
age-dependent branching processes of the Bellman-Harris type, where
individuals have i.i.d. life spans, and split into i.i.d. random
numbers of children, independently of life span, the time structure
thus not being affected by the conditioning.

But what if the flow of time is no longer independent of
reproduction? Even the simplest case, that of a splitting
reproduction at death, but not independently of age at
death/splitting, would seem to offer difficulties, and the same
certainly applies to the more realistic general processes where
reproduction occurs as a point process during life, thus mimicking
the yearly births of wildlife, or the even more erratic reproduction
pattern of humans.

The conceptual framework is intuitive. Starting from an Eve,
individuals live and give birth independently of one another. At
birth each individual inherits a {\em type} from her mother. The
type, in its turn determines the probability measure over all
possible life careers, including a life span and a marked point
process which reports the successive ages at bearing, and the types
of children at the various bearings. Note that multiple births are
not excluded. The branching property can be summarised into the fact
that given her type and birth time, the daughter process of any
individual born is independent of all individuals not in her progeny
(into which she herself is included).

We set out to prove that this branching property remains true for
processes conditioned to die out. Initially, we shall not mention
supercriticality, and only ask that the probability of extinction is
non-zero for any starting type. (If that probability is one, the
conditioning does not change anything!) Largely, the proof is a
matter of conceptual clarity or discipline, which unfortunately
forces us into the somewhat burdensome notation of probabilities on
tree spaces, obscuring the essential simplicity of the matter.

The main idea behind the proof is, however, easily outlined. Indeed,
consider an individual, and condition upon her progeny ultimately
dying out. Her own life career is then affected precisely through
her only being able to have daughters whose progeny in their turn
must ultimately face extinction. In all other respects her life is
independent of all others, once her type is given. This
reestablishes the branching character, but with a suitably amended
probability measure over her life career, which clearly is
non-supercritical in the sense that the probability of ultimate
extinction is one, from any starting type that can be realised.

If the original process is, furthermore, supercritical, {\em i.e.}
has a positive Malthusian parameter, the conditioned process will
turn out to be subcritical, in the sense of the Malthusian parameter
being negative, if it exists.

\section{Notation}

\subsection{The Ulam-Harris family space}

We choose to work within the classical Ulam-Harris framework,
\cite{Jagers:1989}, identifying individuals with sequences of
natural numbers so that $x=(x_1,x_2,\dots,x_n)$ denotes the $x_n$th
child of the \dots\ of the $x_2$th child of the $x_1$th child of the
ancestor. The ancestor is denoted by an ``empty'' sequence $e$
(mnemonic for ``empty'' or ``Eve''), and the set of all possible
individuals is
$$
\T = e\cup\bigcup_{n\in\N}\N^n.
$$
The concatenation of $x,y\in\T$ is   $xy$, and thus $ex=xe=x$ for all $x\in\T$.

For any $e\neq x = (x_1,x_2,\dots,x_n)$ $x$'s \emph{mother} is $mx =
(x_1,\dots,x_{n-1})$, her \emph{rank} in the sibship is $rx = x_n$,
and $x$'s \emph{generation} $g(x)=n$. We agree that $me=re=e$ and
$g(e)=0$. Hence, $mxrx=x$ for $x\in\T$, and $m$ can be iterated so
that $m^nx$ is $x$'s $n$th grandmother, provided $g(x)>n$.

Clearly $x$ \emph{stems from} $y$, usually written $x\succeq y$, if
$m^nx=y$ for some $n\in\N\cup\{0\}$, or equivalently if there exists
a $z\in\T:x=yz$. In this terminology, $x$ stems from herself,
$x\preceq x$. In other words, $(\T,\preceq)$ is a partially ordered
set (a semilattice). We define $x\sim y$ if $x\succeq y$ or
$x\preceq y$, i.e.\ $x$ and $y$ are in direct line of descent.
($\sim$ is not an equivalence relation.)

For $A,B\subseteq\T, x\in\T$, we write $x\succeq A$ if there exists
a $y\in A$ such that $x\succeq y$, and $A\preceq B$ if $x\succeq A$
for all $x\in B$. The \emph{progeny} of $A\subseteq\T$ is defined as
$\mathrm{Pr}\,A=\{x\in\T:x\succeq A\}$.

We call a set $L\subset\T$ a stopping line, or \emph{line} for
short, if no two members of $L$ are in direct line of descent:
$x,y\in L,x\neq y\Rightarrow x\not\sim y$. We say that a line $L$ is
a \emph{covering line} if for all $x\in\T$ there exists a $y\in L$
such that $x\sim y$.

\subsection{Life space and population space}

Let $(\Omega_{\ell},\mathscr{A}_{\ell})$ be a \emph{life space} so
that $\omega\in\Omega_{\ell}$ is a possible life career of
individuals. Any individual property, such as mass at a certain age
or life span, is viewed as a measurable function (with respect to
the $\sigma$-algebra $\mathscr{A}_{\ell}$) on the life space. This
should be rich enough to support, at least, the functions
$\tau(k),\sigma(k)$ for $k\in\N$. Here
$\tau(k):\Omega_{\ell}\to\R_+\cup\{\infty\}$ is the mother's age at
the $k$th child's birth, $0\leq \tau(1)\leq\tau(2)\leq\cdots\leq
\infty$. If $\tau(k)=\infty$, then the $k$th child is never born.
$\sigma(k):\Omega_{\ell}\to\mathcal{S}$ is the child's type,
obtained at birth. The type space $\mathcal{S}$ has a (countably
generated) $\sigma$-algebra $\mathscr{S}$. The whole reproduction
process is then the marked point process $\xi$ with $\xi(A\times
B)=\#\{k:\sigma(k)\in A, \tau(k)\in B\}$ for
$A\in\mathscr{S},B\in\mathscr{B}$, the Borel algebra on $\R_+$.

The population space is defined as $(\Omega,
\mathscr{A})=(S\times\Omega_{\ell}^{\T},\mathscr{S}\times\mathscr{A}_{\ell}^{\T})$.
$U_M$ is the projection $\Omega\to \Omega_{\ell}^M$, for $M\subseteq
\T$. For simplicity $U_x=U_{\{x\}}$ and similarly $\mathrm{Pr}\,x
=\mathrm{Pr}\{x\}$. The following $\sigma$-algebras are important:
$$
\mathscr{F}_L=\mathscr{S}\times\sigma(U_x:x\not\succeq
L)=\mathscr{S}\times\sigma(U_x:x\notin\mathrm{Pr}\,L),
$$
for $L\in\T$. Since $L\preceq M\Rightarrow
\mathrm{Pr}\,L\supseteq\mathrm{Pr}\,M\Rightarrow\mathscr{F}_L\subseteq\mathscr{F}_M$,
it holds that $(\mathscr{F}_L:L\in\T)$ is a filtration under
$\preceq$. In the usual manner, the definition of the
$\sigma$-algebras $\mathscr{F}_L$ can be extended to
$\sigma$-algebras of events preceding random lines $\mathcal{L}$
which are optional in the sense that events $\{\mathcal{L}\preceq
L\} \in \mathcal{F}_L$ \cite{Jagers:1989}.

Functions $\xi,\tau(k)$ and $\sigma(k)$ were defined on the life
space but we want to be able to speak about these quantities
pertaining to a given $x\in\T$. We write $\xi_x=\xi\circ U_x$, $x$'s
reproduction process, $\tau_x = \tau(rx)\circ U_{mx}$, $x$'s
mother's age at $x$'s birth, and $x$'s type
$\sigma_x=\sigma(rx)\circ U_{mx}$. Note the difference between
$\tau(k)$ and $\tau_k$, $\sigma(k)$ and $\sigma_k$, for
$k\in\N\subset\T$.

Finally, the process is anchored in real time by taking Eve to be
born at time 0, and later birth times $t_x,x\in\T$  recursively
determined  by $t_e=0$ and $t_x=t_{mx}+\tau_x$ for $e\neq x\in\T$.
The meaning of $t_x=\infty$ is that $x$ is never born, so that
$\mathscr{R}=\{x\in\T:t_x<\infty\}$ is the set of \emph{realised}
individuals. This set is optional, $\mathcal{F}_{L\cap\mathscr{R}}$
is well defined \cite{Jagers:1989}, and so is the $\sigma$-algebra
$\mathcal{F}_{\mathscr{R}}$ of events pertaining only to realised
individuals. The probability space restricted to such events is that
where a branching processes really lives, {\em cf.} \cite{Neveu},
\cite{Chauvin}.

\subsection{The probability measure and branching property}

The setup is that for each $s\in\mathcal{S}$ there is a probability
measure $P(s,\cdot)$ on the life space
$(\Omega_{\ell},\mathscr{A}_{\ell})$, such that the function $s\to
P(s,A)$ is measurable with $A\in\mathscr{A}_{\ell}$. For any
$s\in\mathcal{S}$ this kernel (the {\em life kernel}) defines a {\em
population  probability} measure $\Prob_s$ on $(\Omega,\mathscr{A})$
with an ancestor of type  $\sigma_e=s$ and such that given
$\sigma_x$, $x$'s life will follow the law $P(\sigma_x,\cdot)$
independently of the rest of the process, see \cite{Jagers:1989}.

Indeed, the basic branching property of the whole process can be
characterised by a generalisation of this in terms of the mappings
$S_x=(\sigma_x,U_{\mathrm{Pr}\,x}):\mathcal{S}\times\T\to
\mathcal{S}\times\T$, which renders $x$ the new Eve. Let
$T_x=S^{-1}_x$ and $\{A_x:x\in L\}\subseteq\mathscr{A}$. Then,
\begin{equation*}
\Prob_s\bigg(\bigcap_{x\in L}T_xA_x\bigg|\mathcal{F}_L\bigg)=
\prod_{x\in L}\Prob_{\sigma_x}(A_x)
\end{equation*}
for lines $L$. This  remains true for optional lines and in particular
\begin{equation}\label{branch_prop}
\Prob_s\bigg(\bigcap_{x\in
L\cap\mathscr{R}}T_xA_x\bigg|\mathcal{F}_{L\cap\mathscr{R}}\bigg)=
\prod_{x\in L\cap\mathscr{R}}\Prob_{\sigma_x}(A_x) ,
\end{equation}
where the intersection over the empty set is taken to be
$\Omega$ and the empty product is ascribed the value one.
The interpretation is that the daughter processes of all realised
individuals $x$ in a line are independent given the prehistory of
the line with the population probability measure $\Prob_{\sigma_x}$,
the only dependence upon the past thus being channelled through the
type $\sigma_x$ and the placing in time $t_x$. This is the branching
property. We shall see that it remains true for processes, bound to
die out.

\section{Conditioning on extinction}

Denote by $E$ the event that the branching process starting from Eve
dies out, i.e. that $\mathscr{R}$ has only a finite number of
elements. Let $q_s=\Prob_s(E)$ and $E_x$ be the event that the
branching process starting from $x$ dies out, $E_x=T_xE$.  Write
$\tilde{\Prob}_s(\cdot)=\Prob_s(\cdot | E)$, which clearly only
makes sense for $s\in \mathcal{S}$ such that $q_s>0$, and let $\tilde{\Ex}_s$
denote expectation with respect to $\tilde{\Prob}_s$.

\begin{theorem}
Any branching process conditioned on extinction remains a branching
process, but with extinction probability one. Its life kernel is
$\tilde{P}(s,A):=\tilde{\Prob}_s(\mathcal{S}\times
A\times\Omega_{\ell}^{\T\setminus\{e\}})$ for
$A\in\mathscr{A}_{\ell}$. Thus, for any covering lines $L$ and
$\{A_x:x\in L\}\subseteq\mathscr{A}$
\begin{equation}\label{cond_branch_prop}
\tilde{\Prob}_s\bigg(\bigcap_{x\in
L\cap\mathscr{R}}T_xA_x\bigg|\mathcal{F}_{L\cap\mathscr{R}}\bigg) =
\prod_{x\in L\cap\mathscr{R}}\tilde{\Prob}_{\sigma_x}(A_x) .
\end{equation}
Furthermore, the Radon-Nikodym derivative $d\tilde{\Prob}_s/d\Prob_s$
with respect to the $\sigma$-algebra $\mathcal{F}_{L\cap\mathscr{R}}$
is given by
\begin{equation}\label{radon-nikodym}
\frac{d\tilde{\Prob}}{d\Prob}\bigg|_{\mathcal{F}_{L\cap\mathscr{R}}}=\frac{1}{q_s}\prod_{x\in L\cap\mathscr{R}}q_{\sigma_x}.
\end{equation}
\end{theorem}

\begin{proof}
First, note that
\begin{equation*}
E= \bigcap_{x\in L\cap\mathscr{R}}T_xE,
\end{equation*}
for covering lines $L$. Indeed, since
$\{L\cap\mathscr{R}=\emptyset\}\subseteq E$ and intersection over an
empty index set yields the full space,
\begin{eqnarray*}\label{E}
  E &=& (E\cap \{L\cap\mathscr{R}=\emptyset\}) \cup (E\cap
\{L\cap\mathscr{R}\neq \emptyset\})\\
  &=& \bigg(\{L\cap\mathscr{R}=\emptyset\}\cap\bigcap_{x\in
L\cap\mathscr{R}}T_xE \bigg) \cup\bigg(
\{L\cap\mathscr{R}\neq\emptyset\}\cap \bigcap_{x\in
L\cap\mathscr{R}}T_xE \bigg)\\
  &=&\bigcap_{x\in L\cap\mathscr{R}}T_xE .
\end{eqnarray*}

The branching property \eqref{branch_prop} implies that
\begin{equation}\label{q_prod}
\Prob_s(E|\mathcal{F}_L)=\prod_{x\in L\cap\mathscr{R}}q_{\sigma_x}
=\Prob_s(E|\mathcal{F}_{L\cap\mathscr{R}}).
\end{equation}
Hence, for any  covering line $L$ and $A\in\mathcal{F}_L$
$$
\tilde{\Prob}_s(A)=\Ex_s\bigg[\Prob_s(E|\mathcal{F}_L);A\bigg]/q_s
=\Ex_s\bigg[\prod_{x\in L\cap\mathscr{R}}q_{\sigma_x};A\bigg]/q_s,
$$
and thus \eqref{radon-nikodym} holds. Equations \eqref{branch_prop} and \eqref{q_prod} yield
\begin{align*}
\Prob_s\bigg(\bigcap_{x\in L\cap\mathscr{R}}T_xA_x\bigg|\mathcal{F}_{L\cap\mathscr{R}}\bigg) &=
\prod_{x\in L\cap\mathscr{R}}q_{\sigma_x}\tilde{\Prob}_{\sigma_x}(A_x) =
\Prob_s(E|\mathcal{F}_{L\cap\mathscr{R}})\prod_{x\in L\cap\mathscr{R}}\tilde{\Prob}_{\sigma_x}(A_x),
\end{align*}
and \eqref{cond_branch_prop} follows.
\end{proof}

\begin{remark} For $L=\mathbb{N}$, the first generation,
$$
\frac{d\tilde{\Prob}_s}{d\Prob_s}\bigg|_{\mathcal{F}_{\N\cap\mathscr{R}}}
=\frac{1}{q_s}\prod_{k\in\mathbb{N}:k\leq X}q_{\sigma_k},
$$
where $X:=\xi_e(\mathcal{S}\times\R_+)$ is Eve's total offspring. In
single-type processes with extinction probability $q$ therefore
\begin{equation}\label{single_rn}
\frac{d\tilde{\Prob}}{d\Prob}\bigg|_{\mathcal{F}_{\N\cap\mathscr{R}}}=q^{X-1}
\end{equation}
and in the Bienaym\'e-Galton-Watson case
$$
\tilde{\Prob}(X=k)=\Ex[q^{X-1};X=k]=\Prob(X=k)q^{k-1},
$$
in perfect agreement with \cite[Theorem I.12.3]{AthreyaNey:1972}.
\end{remark}

\begin{example}
Consider single-type Sevastyanov splitting processes, where
individuals have a life span distribution $G$ and at death split
into $k$ particles with probability $p_k(u)$, if the life span
$L=u$. By \eqref{single_rn} we conclude that
$$
\tilde{\Prob}(X=k, L\in du) =\Ex[q^{X-1};X=k, L\in du] =p_k(u)q^{k-1}G(du).
$$
Hence
$$
\tilde{G}(du)=\tilde{P}(L\in du)=\sum_{k=0}^\infty p_k(u)q^{k-1}G(du)
$$
and $\tilde{\Prob}(X=k, L\in du)=\tilde{p}_k(u)\tilde{G}(du)$, with
$$
\tilde{p}_k(u) = \frac{p_k(u)q^{k-1}}{\sum_{i=0}^{\infty}p_i(u)q^{i-1}},
$$
and the conditioned Sevastyanov process remains a Sevastyanov
process with life span distribution $\tilde{G}$ and splitting
probabilities $\tilde{p}_k(u)$.
\end{example}

Finally, we address the question of whether a supercritical process
conditioned on extinction is subcritical. It is clear that the
conditioned branching process has extinction probability one, for
any starting type, but this would also be the case if the process
were nontrivially critical.

In the single type case it
follows from \eqref{single_rn} that the expected total number of children per
individual in the conditioned process satisfies
$$
\tilde{\Ex}[X]=\Ex[Xq^{X-1}]=f'(q),
$$
in terms of the offspring generating function $f$ of the embedded
Bienaym\'e-Galton-Watson process. It is well known that $f'(q)<1$ if
the original process was supercritical, see \cite{AthreyaNey:1972},
and we conclude that the conditioned process is subcritical.

Also in the general case it is enough to consider the embedded
generation counting process. If $X_n(A) = \mathrm{card}\{x\in\N^n\cap\mathscr{R}: \sigma_x\in A\}$ denotes the number of individuals
in $n$th generation of type $A$ and $q:=\sup_{\mathcal{S}}q_s<1$,
$$\tilde{\Ex}_s[X_n(\mathcal{S})]=\Ex_s\bigg[X_n(\mathcal{S}) e^{\int_{\mathcal{S}}\log q_rX_n(dr)}\bigg]/q_s\le
\Ex_s\bigg[X_n(\mathcal{S})q^{X_n(\mathcal{S})}\bigg]/q_s\to 0,
$$
as $n\to \infty$ (by dominated convergence, since $X_n(\mathcal{S})$ must either tend to zero or to infinity). But the expected size of the embedded process tending to zero means exactly that the process is subcritical.

\end{document}